\documentclass{article}

\usepackage{amsmath,amssymb,latexsym,graphicx}

\newtheorem{theorem}{Theorem}[section]
\newtheorem{lemma}[theorem]{Lemma}

\newcommand{\Nn}[1]{\mathbb{N}^{#1}}  

\newenvironment{proof}{\noindent\textsc{Proof: }}
{\hspace{\stretch{1}}$\Box$\medskip}

\begin{document}

\title
  {A proof of Bj\"orner's totally nonnegative conjecture}

\author{Michael Bj\"orklund \quad\quad
        Alexander Engstr\"om\footnote{{\tt \{mickebj,alexe\}@math.kth.se},
        Department of Mathematics, Royal Institute of Technology, S-100 44 Stockholm, Sweden.}}

\date\today

\maketitle

\abstract{The McMullen Correspondence gives a linear dependence between M-sequences of length $ \lfloor \frac{d}{2} \rfloor + 1 $ 
and $ f $--vectors of simplical $d$--polytopes. Denote the transfer matrix between $ g $ and $ f $ by $ M_d $. Recently, Bj\"orner 
proved that any $ 2 \times 2 $--minor of $ M_d $ is nonnegative and conjectured that the same would be true for arbitrary minors. 
In this paper we answer the question in the affirmative.}

\section{Introduction}

Let $ P $ be a $ d $--dimensional simplicial polytope, and let $ f_i $ denote the 
number of $ i $--dimensional faces. The nonnegative integer vector 
\[
f = ( f_0 , \ldots , f_{d-1} ) 
\]
is called the $ f $-vector of $ P $. We assume the convention $ f_{-1} = f_{d} = 1 $. 
The Euler-Poincar\'e formula gives a first rough estimate on the dependence between the elements 
of $ f $, 
\[
-f_{-1} + f_0 - f_1 + \ldots + (-1)^{d-1} f_{d-1} + (-1)^{d} f_d = 0.
\]
Define the numbers $ g_k = h_{k} - h_{k-1} $ for $ k = 0 , \ldots , \lfloor \frac{d}{2} \rfloor  $, 
where 
\[
h_i = \sum_{j=0}^{i} (-1)^{i+j} \, \binom{d-j}{i-j} \,f_{j-1}, \qquad i = 0, \ldots , d. 
\]
The nonnegative vector $ g =  ( g_0 , \ldots , g_{\lfloor \frac{d}{2} \rfloor} ) $ will be referred
to as the $ g $--vector of the simplical $d$--polytope $ P $. The McMullen correspondence asserts that the map 
\[
g \mapsto g \cdot M_d,
\]
where $ M_d $ is a nonnegative $ ( \lfloor \frac{d}{2} \rfloor + 1 ) \times d $ matrix given by 
\[
m_{ij} = \binom{d+1-i}{d-j} - \binom{i}{d-j} \qquad i = 0, \ldots , \lfloor \frac{d}{2} \rfloor 
\quad \textrm{and} \: j = 0 , \ldots , d-1,
\] 
is a bijection between $ M $--sequences $ g $ with $ g_1 = n - d - 1 $ and $ f $-vectors in 
$ \Nn{d}_0 $ of simplicial $ d $-polytopes with $ n = g_1 + d + 1 $ vertices. Recall that a sequence
$ n_0 , n_1 , \ldots , $ of nonnegative integers is a $ M $-sequence if $ n_0 = 1 $ and 
\[
\partial^k(n_k) \leq n_{k-1} \qquad \textrm{for all } \: k > 1,
\]
where $ \partial^k $ is the ``k-boundary''-operator (see e.g \cite{Z}, 262). A theorem by Macaulay \cite{Ma}
gives an algebro-combinatorial characterization of such sequences. 

This paper is devoted to the  analysis of the transfer matrices $ M_d $ above. In dimensions $ d = 1,2,$ and
$ 3 $, they are given by
\[
M_1 = \left( \begin{array}{cc}
1 & 2
\end{array} \right)
\qquad 
M_2 = \left( \begin{array}{ccc}
1 & 3 & 3 \\
0 & 1 & 1 
\end{array} \right)
\qquad
M_3 = \left( \begin{array}{cccc}
1 & 4 & 6 & 4 \\
0 & 1 & 3 & 2 
\end{array} \right).
\]
Bj\"orner proved in \cite{B} that every $ 2 \times 2 $-minor in $ M_d $ is nonnegative. In the same paper he conjectured 
that any minor in $ M_d $ is nonnegative (such a matrix is called \emph{totally nonnegative}). The conjecture has 
been verified up to dimension $d= 13 $ by A. Hultman. In this paper we settle the conjecture in the affirmative for all dimensions.

\section{Lattice paths and nonnegative minors}

A path from $(x_1,y_1)$ to $(x_2,y_2)$ in
$\mathbb{Z}^2$, where $x_1\leq x_2$ and
$y_1\leq y_2$, is called a lattice path
if only the steps $(1,0)$ and $(0,1)$ are
allowed. The number of lattice paths
from $(0,0)$ to $(m,n)$ which do not
touch the line $y=x+t$ are
${m+n \choose n} - {m+n \choose m-t}$
if $t>0$ \cite{A,M}. This is sometimes
referred to as the ballot numbers.
The weight of a path is the product
of the weights of its arcs.
From any subset $A$
of $\mathbb{Z}^2$ one can construct a
planar acyclic directed graph with vertex set
$A$ and arcs of types $(x,y)\rightarrow
(x+1,y)$ and $(x,y)\rightarrow (x,y+1)$.

\begin{theorem}
 The matrix $M_d$ is
 totally nonnegative for all $d$.
\end{theorem}

\begin{proof}
For integers $n\geq 2$ define the graphs $T_n$
with vertex set
\[\{(x,y)\in\mathbb{Z}^2\mid
  x\leq \lceil n/2 \rceil-1 , y-x\leq \lfloor n/2 \rfloor, \textrm{ and }
  x+y\geq \lceil n/2 \rceil-1 \}.\]
The weight of horizontal arcs in $T_n$ 
is $1$, and the weight of any vertical arc
$(x,y)\rightarrow (x,y+1)$ is $w_y$.
The graphs $T_8$ and $T_9$ are depicted in
figures \ref{fig:p8} and \ref{fig:p9}. 
\begin{figure}
  \begin{center}
  \includegraphics*{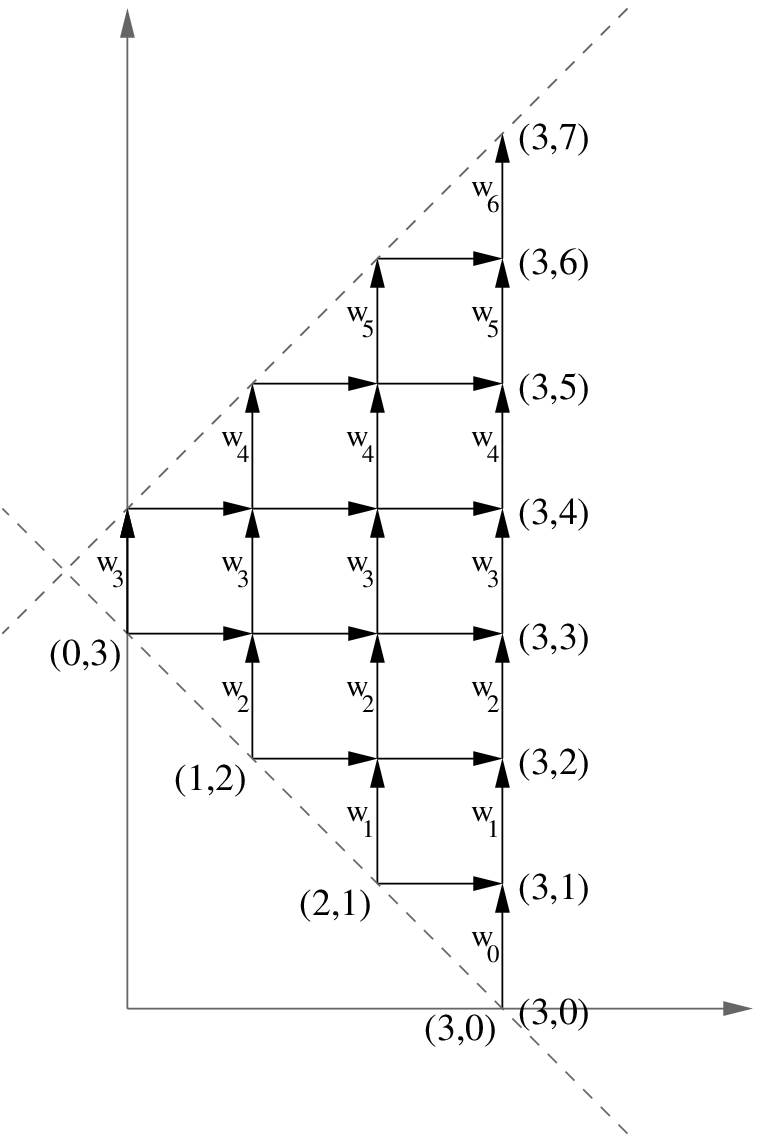}
  \end{center}
  \caption{The weighted planar graph $T_8$}
  \label{fig:p8}
\end{figure}
\begin{figure}
  \begin{center}
  \includegraphics*{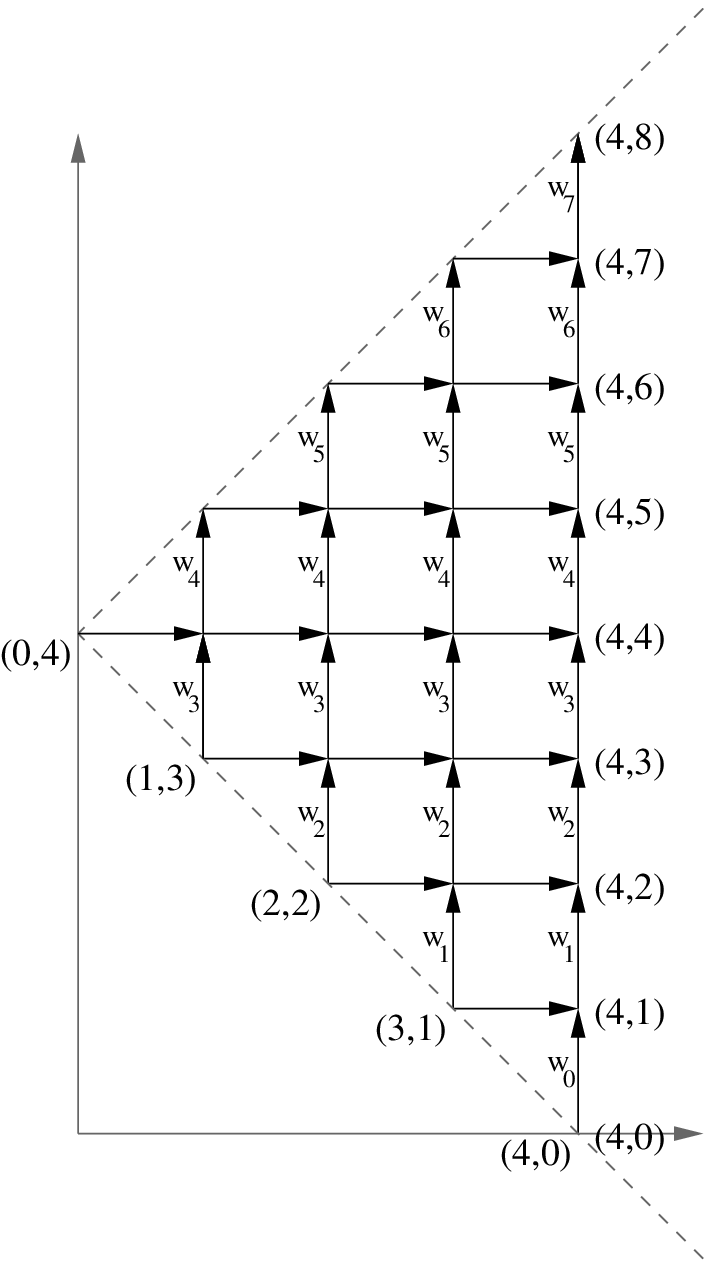}
  \end{center}
  \caption{The weighted planar graph $T_9$}
  \label{fig:p9}
\end{figure}

For all $0\leq i \leq \lceil n/2 \rceil-1$ and
$0\leq j< n$, let the sum of the weights of the
directed paths from $(\lceil n/2 \rceil-1-i,i)$ to
$(\lceil n/2 \rceil-1,j)$ be $W(n,i,j)$. If $i>j$ then
$W(n,i,j)=0$ and otherwise it is
\[w_iw_{i+1}\cdots w_{j-1}\left({j \choose i} - 
{j \choose n-i}\right)\] by the ballot
numbers. 

Now define the weights of the vertical arcs
 as $w_i={n \choose i+1}
/ {n \choose i}$. Note that all arc weights
are positive real numbers. For 
$i\leq j$ we get that
\[\begin{array}{rcl}
 W(n,i,j)  &=& \displaystyle
   \frac{ {n \choose i+1} }{ {n \choose i} }
   \frac{ {n \choose i+2} }{ {n \choose i+1} }
   \cdots
   \frac{ {n \choose j} }{ {n \choose j-1} }
   \left({j \choose i} - 
   {j \choose n-i}\right) \\
& = & \displaystyle
{n \choose j}{n \choose i}^{-1}{j \choose i} -
{n \choose j}{n \choose i}^{-1}{j \choose n-i} \\
& = & \displaystyle
{n-i \choose n-j} - {i \choose n-j}.
\end{array}
\]   
The values of $W(n,i,j)$ is
tabulated in the
$\lceil n/2 \rceil \times n$ matrix $W_n$:
\[ \left(
   \begin{array}{ccccc}
   {n-0 \choose n-0} - {0 \choose n-0}
   &
   {n-0 \choose n-1} - {0 \choose n-1}
   &
   {n-0 \choose n-2} - {0 \choose n-2}
   &
   \cdots
   &
   {n-0 \choose n-(n-1)} - {0 \choose n-(n-1)}\\
   0
   &
   {n-1 \choose n-1} - {1 \choose n-1}
   &
   {n-1 \choose n-2} - {1 \choose n-2}
   &
   \cdots
   &
   {n-1 \choose n-(n-1)} - {1 \choose n-(n-1)}\\
   0
   &
   0 
   &
   {n-2 \choose n-2} - {2 \choose n-2}
   &
   
   &
   {n-2 \choose n-(n-1)} - {2 \choose n-(n-1)}\\
   \vdots & \vdots & && \vdots \\
   0
   &
   0
   &
   0
   &
   &
   {n- \lceil n/2 \rceil +1\choose n-(n-1)} - { \lceil n/2 \rceil -1  \choose n-(n-1)}\\
   \end{array}
   \right)
\]
By removing the leftmost column of $W_n$ we get $M_{n-1}$.
Thus it is sufficient to prove that $W_n$ is totally nonnegative
to conclude the same about $M_{n-1}$. Fomin and Zelevinsky wrote
a nice survey on testing total positivity and related questions
\cite{FM}. We need a result by Lindstr\"om \cite{L}, and
 Gessel and Viennot \cite{GV}.
\begin{lemma} If a weighted acyclic directed planar graph 
has nonnegative real weights, then its weight matrix is
totally nonnegative. 
\end{lemma}
The weight matrix $X$ of a weighted acyclic directed 
planar graph $G$, given a set $I$ of sinks and $J$ of
sources of $G$, is a matrix with the rows indexed by $I$ and
the columns indexed by $J$. On the position $i,j$ of $X$,
where $i\in I$ and $j\in J$, is the sum of the weights of
all paths from the source $i$ to the sink $j$.

The graph $T_n$ with its $\lceil n/2 \rceil$ sources and 
$n$ sinks described earlier gives the weight matrix $W_n$.
All weight are nonnegative, and hence $W_n$ and its
submatrix $M_{n-1}$ are totally nonnegative.
\end{proof}

\bibliographystyle{amsplain}

\end{document}